\documentstyle{article}
\begin{document}
\begin{large}

\centerline{\bf Total subspaces in dual Banach spaces which are not norming}
\centerline{\bf over any infinite dimensional subspace}
\centerline{\bf by}
\centerline{\bf M.I.Ostrovskii}
\centerline{\bf 1. Introduction.}

{\bf Abstract}. The main result: the dual of separable Banach space $X$ 
contains a total subspace which is not norming over any infinite dimensional
subspace of $X$ if and only if $X$ has a nonquasireflexive
quotient space with the strictly singular quotient mapping.

Let $X$ be a Banach space and $X^{*}$ be its dual space.

Let us recall some basic definitions.

A subspace $M$ of $X^{*}$ is said to be {\it total}
if for  every $0\neq x\in X$
there is an $f\in M$ such that $f(x)\neq 0.$

A subspace $M$ of $X^{*}$ is said to be {\it norming over
a subspace} $L\subset X$
if for some $c>0$ we have
$$(\forall x\in L)(\sup_{f\in S(M)}
 |f(x)|\ge c||x||){,}$$
where $S(M)$ is the unit sphere of M. If $L=X$ then $M$ is called
{\it norming}.

The following natural questions arise:

1) How far could total subspaces  be  from  norming  ones?  (Of
course, there are many different concretizations of this question.)

2) What is the structure of Banach spaces, whose duals  contain
total ``very'' nonnorming subspaces?

3) What is the structure of total subspaces?

These questions was studied by many authors: [Al], [B,  p.~208--216],
[BDH], [DJ], [DL], [D], [F], [G], [Ma], [Mc], [M1], [M2],
[O1], [O2], [P], [PP], [S1], [S2].  The  obtained  results  find
applications in the theory of Frechet  spaces  [BDH], [DM], [MM1],
[MM2], [M2]; in the theory of improperly posed problems [O3], [PP,
pp.~185--196] and in the theory of universal bases [Pl, p.~31].

The present paper is devoted to the following natural class  of
subspaces which are far from norming ones. A subspace $M$  of $X^{*}$  is
said to be {\it nowhere norming} if it is not norming over every  infinite
dimensional subspace of $X$. If $X$ is such that $X^{*}$ contains a total
nowhere norming subspace then we shall write $X\in TNNS$. This class  was
introduced by W.J.Davis and W.B.Johnson in [DJ],  where  the  first
example of a total nowhere norming subspace was constructed. In the
same paper it was noted that J.C.Daneman proved that every infinite
dimensional subspace of $l_{1}$ is norming over some infinite dimensional
subspace of $c_{0}$. In [O2] a class of the spaces with $TNNS$ property was
discovered. A.A.Albanese [Al] proved that the spaces  of  type $C(K)$
are not in $TNNS$. The problem of description of  Banach  spaces  with
$TNNS$ property arises in a natural way.

Our main result (Theorem  2.1)  states that for a separable
Banach space $X$ we have $X\in TNNS$  if  and  only  if  for some
nonquasireflexive  Banach  space $Y$  there exists a surjective
strictly singular operator $T:X\to $Y.

Section 3 is devoted to the proof of the auxiliary Theorem
2.4. Using the same method we are able to prove the  following
result (Theorem 3.1):

A Banach  space $M$  is  isomorphic  to  a  total  nonnorming
subspace of the dual of some Banach space if and  only  if $M^{*}$
contains a closed norming subspace of infinite codimension.

Thus the class of total nonnorming subspaces coincide with  the
class of Banach  spaces  that  gives  a  negative  solution  to  the
J.J.Schaffer's problem [Sc, p.~358] (see [DJ, p.~366]).

Section 4 is devoted to several remarks concerning general
(not necessarily separable) spaces, in particular, we  show  that
Banach spaces with the Pelczynski property are not $TNNS$.

Section  5  is   devoted   to   an   example   of   a
nonquasireflexive separable Banach space without  the  Pelczynski
property and without $TNNS$ too.

We hope that our notation is standard and  self-explanatory.
For a  subset $A$  of  a  Banach space $X,\ \hbox{lin}~A,\ A^{\bot}$
and $\hbox{cl}~A$ are, respectively, the linear span of $A$,
the set $\{x^{*}\in X^{*}:(\forall x\in A)(x^{*}(x)=0)\}$
and the closure of $A$ in the strong topology. For a subset $A$ of  a
dual Banach space $X^{*}, w^{*}-\hbox {cl}~A$  and $A^{\top}$  are,
respectively,  the
closure   of  $A$   in   the   weak$^{*}$   topology   and   the    set
$\{x\in X:(\forall x^{*}\in A)(x^{*}(x)=0)\}$. For an operator $T:X\to Y$
the  notation $T|_Z$
denotes the restriction of $T$ to the subspace $Z$ of $X$.

I wish to express gratitude to V.M.Kadets for  his  valuable
advice.

\centerline {\bf 2. Main result.}

Our sources for Banach space basic concepts and results  are
[DS, LT, W]. The unit ball and the unit sphere of Banach space $X$
are denoted by $B(X)$ and $S(X)$ respectively.  The  term  ``operator''
means a bounded linear operator and the term ``subspace''   means  a
closed linear subspace.

Let us recall some definitions.

A Banach space $X$ is called {\it quasireflexive}  if its canonical
image has finite codimension in $X^{**}$.  The  number $\dim(X^{**}/X)$  is
called the {\it order of quasireflexivity}  of $X$  and  is  denoted  by
$\hbox {Ord}~X$.

An operator $T:X\to Y$ is called {\it strictly singular} if the
restriction of $T$ to any infinite dimensional  subspace  of $X$ is
not an isomorphism.

The main result of this paper is the following.

{\bf 2.1.Theorem.} Let $X$ be a separable Banach space.  Then $X\in TNNS$
if and only if for some nonquasireflexive Banach space $Y$ there
exists a surjective strictly singular operator $T:X\to $Y.

Proof. Let us suppose that such an  operator $T$  exists.  We
may assume without loss of generality that $T$ is  a  quotient map.
Then $T^{*}$ is an isometric embedding of $Y^{*}$ into $X^{*}$.  The   subspace
$M_{1}:=T^{*}(Y^{*})$ is  nowhere  norming  because $T$ is  strictly singular.
This space is also not total. Our aim is to find  an isomorphism
$Q:X^{*}\to X^{*}$ which is a small perturbation of the   identity  operator
and is such that under its action $M_{1}$ becomes   a  total  subspace
but remains a nowhere norming one.

The space $Y$ is nonquasireflexive, hence, by [DJ, p.~360] there
exist a weak$^{*}$ null basic sequence $\{y_{n}\}\subset Y^{*}$,
a bounded sequence
$\{g_{n}\}\subset Y^{**}$ and a partition $\{I_{n}\}^{\infty }_{n=1}$
of  the integers into pairwise disjoint infinite subsets such that
$$g_{k}(y_{n})=\cases{1,&if $n\in I_k$;\cr
0, &if $n\not\in I_{k}$.\cr}$$

Let us denote the vector $T^{*}y_{n}$ by $u^{*}_{n}$.
The operator $T^{*}$ is weak$^{*}$
continuous and isometric  hence $\{u^{*}_{n}\}^{\infty }_{n=1}$
is a weak$^{*}$ null basic
sequence in $X^{*}$ and in $X^{**}$
there exists a bounded sequence $\{v^{**}_{n}\}^{\infty }_{n=1}$
such that
$$v^{**}_{k}(u^{*}_{n})=\cases{1, &if $n\in I_{k}$;\cr
0, &if $n\not\in I_{k}$.\cr}$$

Let $\{s^{*}_{k}\}^{\infty }_{k=1}$ be a normalized sequence
spanning a total subspace in $X^{*}$. Let the operator
$Q:X^{*}\to X^{*}$ be given by
$$
Q(x^{*})=x^{*}+\sum^{\infty }_{k=1}4^{-k}
v^{**}_{k}(x^{*})s^{*}_{k}/||v^{**}_{k}||.
$$

It is clear that $Q$ is an isomorphism. Let $M=Q(M_{1})$. We shall
show that $M$ is a total nowhere norming subspace. Let $0\neq x\in X$ and let
$k\in {\bf N}$ be such that $s^{*}_{k}(x)\neq 0.$
Since $\{u^{*}_{n}\}^{\infty }_{n=1}$ is weak$^{*}$ null then we can
choose $n\in I_{k}$ such that
$|u^{*}_{n}(x)|<4^{-k}s^{*}_{k}(x)/||v^{**}_{k}||$. We have
$$(Q(u^{*}_{n}))(x)=
u^{*}_{n}(x)+4^{-k}s^{*}_{k}(x)/||v^{**}_{k}|| \neq 0. $$
Hence $M$ is total.

Recall that if $U,V$ are subspaces  of a Banach  space,  then
the number
$$
\delta (U,V)=\inf\{||u-v||:u\in S(U),v\in V\}
$$
is called the {\it inclination} of $U$ to $V$.

We shall prove that $M\subset X^{*}$ is nowhere norming. Let us  suppose
that it is not the case and let an infinite dimensional  subspace
$L\subset X$ be such  that $M$  is norming over $L$. By strict  singularity
of $T$  there  is  no any infinite dimensional subspace of $X$
with nonzero inclination  to $\ker T$.  Using  standard  reasoning
with basic sequences (see [Gu]) we can find in $L$ a normalized basic
sequence $\{z_{i}\}$ such that for some sequence $\{t_{i}\}\subset \ker T$
we shall
have $||z_{i}-t_{i}||\le 2^{-i}$, and moreover we may require that
$$
(\forall n\in {\bf N})(\lim_{i\to \infty }s^{*}_{n}(t_{i})=0).\eqno(1)
$$
Let $c>0$ be such that

$$
(\forall x\in L)(\exists f\in S(M))(|f(x)|\ge c||x||).
$$

In particular
$$
(\forall i\in {\bf N})(\exists f_{i}\in S(M))(|f_{i}(z_{i})|\ge c).
$$

By the definition of $M$ we can find such $y^{*}_{i}\in Y^{*}$ that
$$
f_{i}=T^{*}y^{*}_{i}+\sum^{\infty }_{k=1}4^{-k}v^{**}_{k}
(T^{*}y^{*}_{i})s^{*}_{k}/||v^{**}_{k}||.
$$

>From this equality we obtain
$$
||f_{i}||\ge(2/3)||T^{*}y^{*}_{i}||.
$$
Using this inequality we obtain for every  positive  integer
$i$ that
$$
c\le |f_{i}(z_{i})|\le |f_{i}(z_{i}-t_{i})|+|f_{i}(t_{i})|\le
$$
$$
2^{-i}+
|\sum^{\infty }_{k=1}4^{-k}v^{**}_{k}(T^{*}y^{*}_{i})s^{*}_{k}(t_{i})/
||v^{**}_{k}|| |\le
$$
$$
2^{-i}+(3/2)\sum^{\infty }_{k=1}4^{-k}|s^{*}_{k}(t_{i})|.
$$
Using (1) and the boundedness of sequences $\{s^{*}_{n}\}$ and $\{t_{i}\}$  we
arrived at a contradiction. Hence $M$ is nowhere norming.

Now we begin to prove  the  converse  statement.  We  need  the
following result, which easily follows from the arguments of [DJ, p.~358].

{\bf 2.2.Lemma.} Let $X$ be a separable Banach space and let $N$  be  a
subspace of $X^{*}$ such that the strong closure of the  canonical  image
of $X$ in $N^{*}$ is of infinite codimension.
Then $N$ contains a weak$^{*}$ null
basic sequence $\{u^{*}_{n}\}^{\infty }_{n=1}$ such that for
some bounded sequence $\{v^{**}_{k}\}^{\infty }_{k=1}$
in $X^{**}$ and some partition
$\{I_{k}\}^{\infty }_{k=1}$ of the positive integers into
pairwise disjoint infinite subsets we have
$$
v^{**}_{k}(u^{*}_{n})=\cases{1, &if $n\in I_{k}$;
\cr 0, &if $n\not\in I_{k}$.\cr}
$$

It turns out that a total nowhere norming subspace need not  satisfy
the condition of Lemma 2.2.

{\bf 2.3.Proposition.} There exists a total nowhere norming subspace
$L$ of $(l_{1})^{*}$ such that the canonical
image of $l_{1}$ is dense in $L^{*}$.

Proof. Let a Banach space $X$  be  such  that $X^{*}$  is  separable,
contains closed norming subspaces of infinite codimension  and  does
not contain subspaces isomorphic to $l_{1}$. We may take  e.\ g.\
$X=(\sum \oplus J)_{2}$, where $J$ is James' space [LT, p.~25].

We need the following definition. Let $a\ge 0, b\ge 0.$ We shall say
that subset $A\subset X^{*}$ is $(a, b)$-{\it norming} if the following
conditions are satisfied:
$$
(\forall x\in X)(\sup \{|x^{*}(x)|:x^{*}\in A\}\ge a||x||);
$$
$$
\sup \{||x^{*}||:x^{*}\in A\}\le b.
$$

Let $K\subset X^{*}$ be a closed norming subspace of infinite  codimension.
Let $\alpha :l_{1}\to K$  be  some  quotient  mapping.
Hence  the  set $\alpha (B(l_{1}))$  is
$(c,1)$-norming for some $c>0.$ Let a
sequence $\{z_{i}\}^{\infty }_{i=1}\subset X^{*}$ be  such  that
its image  under  the  quotient  map $X^{*}\to X^{*}/K$  is  minimal.  Let  us
introduce the operator $\beta :l_{1}\to X^{*}$ in the following way:
$$
\beta (\{a_{i}\}^{\infty }_{i=1})=(c/2)\sum^{\infty }_{i=1}a_{i}z_{i}/||z_{i}||+
\alpha (\{a_{i}\}^{\infty }_{i=1}).
$$
It is clear that $\beta $  is  injective
and  that  the  set $\beta (B(l_{1}))$  is
$(c/2,1+c/2)$-norming.
Let (finite or infinite) sequence $\{y_{i}\}^{k}_{i=1}\subset X^{*}$ be
such that its image under the quotient
mapping $X^{*}\to X^{*}/\hbox {cl}\beta (l_{1})$ is  minimal
and that
$X^{*}=\hbox {cl}(\hbox {lin}(\{y_{i}\}^{k}_{i=1}\cup \beta (l_1)))$.
Let us represent $l_{1}$ as $l_{1}\oplus l^{k}_{1}$ (or
$l_{1}\oplus l_{1}$ if $k$ equals infinity) and define the
operator $\gamma :l_{1}\oplus l^{k}_{1}\to X^{*}$  in
the following way:
$$
\gamma (\{a_{i}\}^{\infty }_{i=1},\{b_{i}\}^{k}_{i=1})=
\beta (\{a_{i}\}^{\infty }_{i=1})+\sum^{k}_{i=1}b_{i}y_{i}/||y_{i}||.
$$
It is clear that $\gamma $ is injective,
its image is dense in $X^{*}$  and  that
$\gamma (B(l_{1}))$ is $(c/2,1+c/2)$-norming.  Besides  this $\gamma$
is  a  strictly singular operator
since $X^{*}$ does not contain subspaces isomorphic  to
$l_{1}$.

Let $L=\gamma ^{*}(X)\subset (l_{1})^{*}$.
This subspace is total since $\gamma $ is injective.
Since $\gamma (B(l_{1}))$ is $(c/2,1+c/2)$-norming
then $\gamma^{*}|_{X}$  is  an  isomorphic
embedding. Therefore the strict singularity of $\gamma $ implies that $L$  is
nowhere norming. On the other hand it is  easy  to  check  that  the
canonical image of $l_{1}$  in $L^{*}$
may  be  identified  with $\gamma (l_{1})$  and
therefore is dense. The proposition is proved.

In order to make Lemma 2.2 applicable for our purposes we  need
the following result.

{\bf 2.4.Theorem.} Let $X$ be a Banach  space  and $M$  be  a  total
nowhere norming subspace of $X^{*}$. Then there  exists  an  isomorphic
embedding $E:M\to X^{*}$ such that $E(M)$ is a also nowhere norming subspace
and the closure of $E^{*}(X)$ in the strong topology has  infinite
codimension in $M^{*}$.

We postpone the proof until section 3.

Let $X\in TNNS$ and let $M\subset X^{*}$ be a total
nowhere  norming  subspace.
Applying Theorem 2.4 we find an embedding $E:M\to X^{*}$ such that $N=E(M)$ is
a nowhere norming subspace satisfying the condition  of  Lemma  2.2.
Let $\{u^{*}_{n}\}^{\infty }_{n=1}\subset E(M),
\{v^{**}_{k}\}^{\infty }_{k=1}\subset X^{**}$
and $\{I_{k}\}^{\infty }_{k=1}$ be sequences obtained  by
application of Lemma 2.2 to $N=E(M)$.

We need the following definition [JR].

A sequence $\{x^{*}_{n}\}^{\infty }_{n=1}\subset X^{*}$
is called {\it weak}$^{*}$ {\it basic} provided that there
is a sequence $\{x_{n}\}^{\infty }_{n=1}\subset X$
so that $\{x_{n},x^{*}_{n}\}$ is biorthogonal and for each
$x^{*}\in w^{*}-\hbox {cl}(\hbox {lin}\{x^{*}_{n}\}^{\infty }_{n=1})$,
$$
x^{*}=w^{*}-\lim_{n\to \infty }\sum^{n}_{i=1}x^{*}(x_{i})x^{*}_{i}.
$$
By [JR, p.~82] (see also [LT,p.~11]) every bounded away from  0
weak$^{*}$ null sequence in the dual of a separable Banach  space  has  a
weak$^{*}$ basic subsequence. Therefore, we  may  select  a  weak$^{*}$  basic
subsequence
$\{u^{*}_{n(j)}\}^{\infty }_{j=1}\subset \{u^{*}_{n}\}^{\infty }_{n=1}$.
Moreover, by arguments  of [JR]  we
may suppose that the intersection $I_{k}\cap \{n(j)\}^{\infty }_{j=1}$
is infinite for every $k\in {\bf N}$.

Let $Y=X/((\{u_{n(j)}\}^{\infty }_{j=1})^{\top})$
and let $T:X\to Y$ be the quotient map.  The
dual $Y^{*}$ may be naturally identified with
$w^{*}-\hbox {cl}(\hbox {lin}\{u^{*}_{n(j)}\}^{\infty }_{j=1})$.  The
space $Y$  is  nonquasireflexive  because  intersections
$I_{k}\cap \{n(j)\}^{\infty }_{j=1} (k\in {\bf N})$
are infinite. By the  well-known  properties  of  weak$^{*}$  basic
sequences [LT, p. 11] it follows that for some $\lambda <\infty $ we have
$$
B(Y^{*})\subset \lambda w^{*}-\hbox{cl}
(B(\hbox{lin}\{u^{*}_{n(j)}\}^{\infty }_{j=1}))\subset
\lambda w^{*}-\hbox{cl}(B(E(M))).
$$
Therefore $Y^{*}$ is nowhere norming subspace of $X^{*}$.
Hence $T$ is  strictly
singular. The proof of Theorem 2.1 is complete.

{\bf 2.5. Corollary.} If $X$ is a separable Banach space which contains a
complemented subspace $Y$ with $Y\in TNNS$ then $X\in TNNS$.

Proof. The composition of a strictly singular surjection $T:Y\to Z$
and any projection $P:X\to Y$ is a required surjection.

{\bf 2.6. Corollary.} For  every  separable  Banach  space $X$  we have
$X\oplus l_{1}\in TNNS$.

To prove this we  need  only  to  recall  that  there  is  a
surjective strictly singular operator $T:l_{1}\to c_{0}$ [LT, p.~75,~108].

{\bf 2.7. Remark.} There exists a space $X$ for which $X\in TNNS$ but $X$
does not contain any subspaces isomorphic to $l_{1}$. In fact, by the
James-Lindenstrauss theorem [LT, p.~26] there exists  a  separable
space $Z$  for  which $Z^{**}/Z$  is  isomorphic  to $c_{0}$  and $Z^{***}$  is
isomorphic to $Z^{*}\oplus l_{1}$.
Let $X=Z^{**}$. Then there exists a quotient  map
$T:X\to c_{0}$. This map must be strictly singular because if not then $X$
must contain a subspace isomorphic to $c_{0}$ [LT,  p.~53]  but  this
contradicts the fact that $X$ is a separable dual space [LT, p.~103].
At the same time $X$ does not contain subspaces isomorphic to
$l_{1}$ because $X^{*}$ is separable.

\centerline {\bf 3. Total nonnorming subspaces in dual Banach spaces.}

In this section we shall prove Theorem 2.4 and the following
characterization of total nonnorming subspaces.

{\bf 3.1.Theorem.} A Banach space $M$  is  isomorphic  to  a  total
nonnorming subspace of the dual of some Banach  space  if and
only if $M^{*}$  contains  a  closed  norming  subspace  of  infinite
codimension.

We need the following lemmas.

{\bf 3.2.Lemma} [B, p.~39]. If $U$ and $V$ are Banach spaces and
$P:U\to V$ is an operator with nonclosed image  then  the  closure  of
$P(B(U))$ in the strong topology does not contain interior points.

{\bf 3.3.Lemma} [LT,  p.~79].  If $P:U\to V$  is  an  operator  with
nonclosed image and $F:U\to V$ is a  finite  rank  operator  then  the
operator $(P+F)$ has nonclosed image.

{\bf 3.4.Lemma.} Let $P:U\to V$ be an operator  with  nonclosed  image
and let $\varepsilon >0.$ Then there exist a
functional $f\in V^{*}$ and an operator
$P_{1}:U\to V$ such that $f$ does not vanish on ${\rm im}~P$,
the image of $(P-P_{1})$
is one-dimensional, $||P-P_{1}||\le\epsilon $ and
${\rm im}~P_{1}\subset\ker f\cap {\rm cl}({\rm im}~P)$.

Proof. By Lemma 3.2 the closed convex set ${\rm cl}(P(B(U)))$  does
not have interior points in the subspace
$V_{0}={\rm cl}({\rm im} P)\subset V$. Therefore
there exists a functional $f_{0}\in S(V^{*}_{0})$ such that
$$
(\forall v\in P(B(U)))(|f_{0}(v)|\le \varepsilon /2).
$$
It is clear that $f_{0}$ does not vanish on ${\rm im}~P$.
Let $v_{0}\in V_{0}$  be
such that $f_{0}(v_{0})=1$ and $||v_{0}||\le 2.$
Let us define an operator $P_{1}:U\to V$ by
$P_{1}(u)=P(u)-f_{0}(P(u))v_{0}$.
Let $f$ be any continuous  extension  of $f_{0}$
onto the whole $V$. It can be  directly  verified  that $P_{1}$  and $f$
defined in such a way satisfy all the requirements of Lemma 3.4.

Let $X$ be a Banach space and let $M$ be a subspace of its dual.
Every element of $X$ may be considered as a  functional  on  $M$.  So
there is a natural map of $X$ into $M^{*}$. We shall denote this map  by
$H$.

{\bf 3.5.Proposition.} Let $X$ be a Banach space and  let $M$  be  a
total nonnorming subspace in $X^{*}$. Then there exists an  isomorphic
embedding $E:M\to X^{*}$ such that the closure of the $E^{*}(X)$ in the strong
topology is of infinite codimension in $M^{*}$ and the difference
$(E^{*}|_{X}-H)$ is a nuclear operator.

Proof. In our case the map $H$ is injective because $M$ is total
and is not an isomorphic embedding because $M$  is  nonnorming.  By
the open mapping theorem the image of $H$ is nonclosed.

Let us apply Lemma 3.4 to $P=H$  and $\varepsilon =1/4,$  and  denote  the
obtained functional by $f_{1}$ and the obtained operator by $H_{1}$. By Lemma
3.3  the  operator $H_{1}$  also  has  nonclosed  image. Applying Lemma
3.4 to $P=H_{1}$ and $\varepsilon =1/8$
we find functional $f_{2}$ and operator $H_{2}$. We
continue in an obvious way.

We have $||H_{i-1}-H_{i}||<2^{-i-1}$.
Therefore the sequence $\{H_{i}\}^{\infty }_{i=1}$ is
uniformly convergent. Let us denote by $R$ its limit.

The operator $(R-H)$ is nuclear and satisfies the inequality
$$
||R-H||<2^{-1}\eqno(2)
$$
It is clear that the set $H(B(X))$ is (1,1)-norming.  By  this
and  by  inequality  (2)  we  obtain  that  the  set $R(B(X))$  is
(1/2,3/2)-norming. Moreover we have
${\rm cl}({\rm im} R)\subset \bigcap^{\infty }_{i=1}\ker f_{i}$.  The
sequence $\{f_{i}\}^{\infty }_{i=1}$  is  linearly   independent
because   it   is constructed in such a way
that $f_{i+1}$ does not  vanish  on
$\bigcap^{i}_{k=1}\ker f_{k}$.
Therefore ${\rm cl}({\rm im} R)$ is of infinite codimension in $M^{*}$.
Let us introduce an operator $E:M\to X^{*}$ by $(E(m))(x)=(R(x))(m)$.
This operator is an isomorphic embedding because the set $R(B(X))$
is (1/2,3/2)-norming.

It is easy  to  see  that  the  restriction  of $E^{*}$  onto $X$
coincides with $R$. Therefore $(E^{*}|_{X}-H)$ is nuclear
and ${\rm cl}(E^{*}(X))$ is of infinite codimension in $M^{*}$.
The proof is complete.

Proof of Theorem 2.4. A nowhere norming subspace is of course
nonnorming. So we can apply Proposition 3.5. It should  be  noted
that for nowhere norming $M$ the operator $H$ is strictly singular.

Let $E:M\to X^{*}$ be the operator constructed in  Proposition  3.5.
We need only to check that $E(M)$  is  nowhere  norming.  But  this
follows immediately from the fact that $E^{*}|_{X}=R=H+(R-H)$ is strictly
singular as a sum of two strictly singular operators [LT, p.~76].

Proof of Theorem 3.1. The necessity follows immediately from
Proposition 3.5 and the fact that ${\rm cl}(E^*(X))$ is norming subspace
in $M^{*}$.

Let us suppose that $M$ is a  Banach  space  for  which  there
exists a closed norming subspace $V\subset M^{*}$ of infinite codimension. Let
$\{z_{i}\}^{\infty }_{i=1}$be a normalized basic sequence
in $M^{*}/V$ and let $m^{*}_{i}\in M^{*} (i\in {\bf N})$ be
such that $||m^{*}_{i}||\le 2$ and $Q(m^{*}_{i})=z_{i}$,
where $Q:M^{*}\to M^{*}/V$ is a quotient map.

Let $X=V\oplus l_{1}$. Let us define an operator $H:X\to M^{*}$ by
$$
H(v,\{a_{i}\}^{\infty }_{i=1})=v+\sum^{\infty }_{i=1}(a_{i}/i)m^{*}_{i}.
$$

It is clear that this operator is injective but  is  not  an
isomorphic embedding.

The restriction of $H^{*}$ onto $M$  is  an  isomorphic  embedding
because $V$ is a norming subspace. The subspace
$H^{*}(M)\subset X^{*}$  is  total
because $H$ is injective and is not norming because $H$  is  not  an
isomorphic embedding. This  completes the proof of Theorem 3.1.

{\bf 3.6. Corollary.} If $M$ is a total subspace  of $X^{*}$  and $M$  is
quasireflexive then $X$ is quasireflexive and ${\rm Ord}(X)={\rm Ord}(M)$.

Proof. It is known [CY] that  ${\rm Ord}(X^{*})={\rm Ord}(X)$  and  that  the
order of quasireflexivity of a subspace is not greater  than  the
order of quasireflexivity of the whole space.

It is well-known and is  easy  to  see  that  the  duals  of
quasireflexive spaces do not contain norming subspaces of infinite
codimension. Therefore by Theorem 3.1 the subspace $M\subset X^{*}$ is norming.
Hence $X$  is  isomorphic  to  a  subspace  of $M^{*}$.   Using   the
abovementioned result we
obtain ${\rm Ord}(X)\le {\rm Ord}(M^{*})={\rm Ord}(M)$.

Using the abovementioned result once more we
obtain ${\rm Ord}(M)\le {\rm Ord}(X^{*})={\rm Ord}(X)$. The proof is complete.

{\bf 3.7.Remark.} By [DJ, p.~355] nonquasireflexivity of $X$ does  not
yield the existence in $X^{*}$ of an  infinite  codimensional  norming
subspace. Therefore there exist  nonquasireflexive  spaces  which
are not isomorphic to total nonnorming subspaces.

\centerline {\bf 4. Remarks on the nonseparable case and on spaces with}
\centerline {\bf the Pelczynski property.}

Theorem 2.1 and Corollaries 2.5 and 2.6 are not valid in
nonseparable case. In order to prove this let us  show  that  the
space $X=l_{1}\oplus l_{2}(\Gamma )$ does not have $TNNS$
property if ${\rm card}(\Gamma )>2^{c}$.

Let $M$ be a total subspace in $X^{*}$.
Then $2^{\hbox{card}(M)}\ge {\rm card}(X)\ge {\rm card}(\Gamma )$.
Consequently ${\rm card}(M)>c$. Therefore $M$ contains a set of functionals of
cardinality greater than $c$, whose restrictions to $l_{1}$ coincide.
Therefore the intersection of $M$ with the subspace of $X^{*}$ which
vanishes on $l_{1}$ is an infinite dimensional
subspace in $\{0\}\oplus l_{2}(\Gamma )$. If
we ``transfer'' this subspace into $X$ then we shall obtain a subspace
over which $M$ is norming.

{\bf Problem.} Characterize $TNNS$ in the nonseparable setting.

At the moment it is known [Al] that $C(K)\not\in TNNS$ for every compact $K$.

{\bf 4.1.Proposition.} Let Banach space $X$  be  such  that  every
strictly singular operator $T:X\to Y$ is weakly compact. Then $X\not\in TNNS$.

Proof. Evidently it is sufficient to consider the case  when
$X$ is nonreflexive. Let $M$ be a total nowhere norming  subspace  in
$X^*$. Let $X_{M}$ be the completion of $X$ under the norm
$$
||x||_{M}=\sup \{|f(x)|:f\in S(M)\}.
$$
Let $T:X\to X_M$ be the natural embedding. The operator $T$  is
strictly singular because $M$ is nowhere norming. Hence $T$ is weakly
compact.

The subspace $M\subset X^{*}$  may be considered also as a subspace of
$(X_{M})^{*}$. Moreover the
restriction of $T^{*}:(X_{M})^{*}\to X^{*}$ onto $M$ is an isometry.

The operator $T^{*}$ is  weakly  compact  by  V.Gantmacher's  theorem
[DS, VI.4.8]. Therefore $M$ is reflexive, hence the subspace $M$ is weak$^{*}$
closed in $X^{*}$ by the M.~Krein--V.~Smulian theorem [DS, V.5.7]. Since $M$
is a total subspace of $X^{*}$ we obtain $M=X^{*}$. This contradiction
completes the  proof.

{\bf 4.2. Remark.} For separable spaces this proposition  follows
immediately from Theorem 2.1.

The conditions of Proposition 4.1 are satisfied by spaces with
the Pelczynski property. Let us recall the definition.

A Banach space $X$ has the {\it Pelczynski property} if for every
subset $K\subset X^{*}$  that is  not relatively weakly compact there exists a
weakly  unconditionally convergent series $\sum^{\infty }_{n=1}x_{n}$
in $X$ such that
$$\inf _{n}\sup _{x^*\in K}x^{*}(x_{n})>0.$$

This property was introduced by A.Pelczynski in [Pe] (under the
name ``property (V)''). In the same paper it was proved that for any
compact Hausdorff space $S$ the space $C(S)$ has property $(V)$. For other
spaces with the Pelczynski property see [W, pp.~166-172].

The  fact  that  spaces  with the Pelczynski property satisfy
the conditions of Proposition 4.1 follows from the next proposition.

{\bf 4.3.Proposition} [W,p.~172]. Suppose $X$ has the Pelczynski
property. Then for every  operator $T:X\to Y$  that  is  not  weakly
compact there exists a subspace $X_{1}\subset X$ such
that $X_{1}$  is  isomorphic
to $c_{0}$ and the restriction of $T$ onto $X_{1}$ is an isomorphic embedding.

The spaces  which  have  no  $TNNS$  property  need  not  have
the Pelczynski property and need not satisfy the  conditions  of
Proposition 4.1.

A corresponding example is given by James` space $J$. Let  us
recall its definition [LT, p.~25]. The space $J$  consists  of  all
sequences of scalars $x=(a_{1}, a_{2},\ldots , a_{n},\ldots )$ for which
$$
||x||=\sup 2^{-1/2}((a_{p_1} -a_{p_2})^2+(a_{p_2} -a_{p_3})^2+\ldots
+(a_{p_{m-1}} -a_{p_m})^2)^{1/2}<\infty
$$
where the supremum is taken over all choices of $m$ and
$p_1<p_2<\ldots <p_m$; and
$$
\lim_{n\to\infty}a_{n}=0.
$$

It is easy to see that the operator $T:J\to c_{0}$ which maps  every
sequence from $J$ onto the same sequence  in $c_{0}$  is a non-weakly
compact strictly singular operator. On the other hand  any  total
subspace in $J^{*}$ is norming over $J$. This  follows  from  well-known
properties of quasireflexive spaces.

In fact there are nonquasireflexive spaces of this  type  as
is shown in the next section.

\centerline {\bf 5. A nonquasireflexive separable Banach space without the}
\centerline {\bf Pelczynski property whose dual does not contain nowhere}
\centerline {\bf norming subspaces.}

{\bf 5.1.Theorem.} There exists a nonquasireflexive Banach space
$X\not\in TNNS$ and such that there exists a strictly singular  non-weakly
compact operator $T:X\to c_{0}$.

Proof. Let $X=(\sum^{\infty }_{n=1}\oplus J)_{2}$.
The  unit  vectors  in $J$  we  shall
denote by $\{e_{i}\}^{\infty }_{i=1}$. It is known [LT, p.~25]
that $\{e_{i}\}$ is shrinking
basis of $J$, therefore its biorthogonal
functionals $\{e^{*}_{i}\}^{\infty }_{i=1}$ form a
basis of $J^{*}$.

It is clear that vectors
$$
e_{n,j}=(0,\ldots ,0,e_j,0,\ldots ).
$$
(where $e_j$ is on the $n$-th place) after
any numeration preserving order in sequences $\{e_{n,j}\}^{\infty }_{j=1}$ form
a basis of X. We need the following two lemmas about $X$ and its dual.

{\bf 5.2.Lemma.} Every weakly null
sequence $\{x_{m}\}^{\infty }_{m=1}$ in $X$ for which
$\inf||x_{m}||>0$ contains a subsequence equivalent to  the  unit  vector
basis of $l_{2}$.

{\bf 5.3.Lemma.} Every infinite dimensional subspace of $X^{*}$ contains
a subspace isomorphic to $l_{2}$.

This lemmas easily follows by well-known arguments (see [An] and [HW]).

Let us consider an
operator $T:(\sum \oplus J)_{2}\to c_{0}=(\sum \oplus c_{0})_{0}$  defined  by
$$T(x_1,\ldots ,x_n,\ldots )=(x'_1,\ldots ,x'_n,\ldots ),$$
where $(x_{i})$  is a  sequence of
elements of $J$ and $(x'_i)$ is a sequence of elements of $c_{0}$  with  the
same coordinates. It is clear that $T$ is a continuous operator. It
is  not  weakly  compact  because  for  any $n\in {\bf N}$   the   sequence
$(T(\sum^{k}_{j=1}e_{n,j}))^{\infty }_{k=1}$ does
not have limit points in weak topology.  At
the same time the operator $T$ is strictly singular  because  by
Lemma 5.2 the space $X$ does not contain subspaces isomorphic to $c_{0}$.

Let us suppose that $X\in TNNS$. Then by Theorem 2.1 there exists a
surjective strictly singular operator $T:X\to Z$ where $Z$ is a certain
Banach space. Consequently $Z^{*}$ is isomorphic to a subspace of $X^{*}$.
By Lemma 5.3 the space $Z^{*}$ contains a subspace isomorphic to $l_{2}$.  Let
us denote this subspace by $U$. Let $R$ be the
quotient  map $R:Z\to Z/U^{\top }$.
The space $(Z/U^{\top })^{*}$ may be in natural way
identified  with $w^{*}-{\rm cl}U$.
Since $U$ is reflexive then by M.Krein-V.Smulian theorem [DS,  V.5.7]
we have $w^{*}-{\rm cl}(U)=U$.
Therefore $Z/U^{\top }$ is isomorphic to $l_{2}$.
Let $\{u_{i}\}^{\infty }_{i=1}$
be sequence in $Z/U^{\top }$ equivalent to
the unit vector basis of $l_{2}$. By
Lemma 2 of [GR] we can find in $X$ a weakly
null sequence $\{x_{i}\}^{\infty }_{i=1}$ for
which $\{RTx_{i}\}^{\infty }_{i=1}$ is a
subsequence of $\{u_{n}\}$. By Lemma 5.2 the sequence
$\{x_{i}\}$ contains a subsequence $\{x_{n_i} \}^{\infty }_{i=1}$
which is equivalent to the unit
vector basis of $l_{2}$. The restriction of $RT$ onto the closed linear
span of $\{x_{n_i}\}^{\infty }_{i=1}$ is an isomorphism.
Because $T$ is strictly singular
this gives us a contradiction.

\centerline {REFERENCES}

[Al] A.A.Albanese, On total subspaces in duals  of  spaces  of  type
$C(K)$ or $L^{1}$, preprint.

[An] A.Andrew,James' quasi-reflexive space is  not  isomorphic  to
any subspace of its dual, Israel  J. Math. 38(1981), 276--282.

[B]  S.Banach, Theorie des operations  lineaires,   Monografje
Matematyczne. No. 1 (Warszawa, l932).

[BDH] E.Behrends, S.Dierolf and P.Harmand, On a problem of  Bellenot
and Dubinsky, Math. Ann. 275 (1986), 337--339.

[CY] P.Civin and  B.Yood,  Quasi-reflexive  spaces,  Proc. Amer.
Math. Soc. 8 (1957), 906--911.

[DJ] W.J.Davis  and  W.B.Johnson, Basic  sequences   and   norming
subspaces in non-quasi-reflexive Banach spaces , Israel J.  Math.
14 (1973), 353--367.

[DL] W.J.Davis and J.Lindenstrauss, On total  nonnorming  subspaces,
Proc. Amer. Math. Soc. 31 (1972), 109--111.

[DM] S.Dierolf and V.B.Moscatelli, A  note  on  quojections,  Funct.
Approx. Comment. Math. 17 (1987), 131--138.

[D]  J.Dixmier, Sur un theoreme de Banach, Duke Math. J. 15 (1948),
1057--1071.

[DS] N.Dunford and J.T.Schwartz, Linear operators. Part I: General
theory (Interscience Publishers, New York - London, 1958).

[F]  R.J.Fleming, Weak$^{*}$-sequential closures and  the  characteristic
of subspaces  of  conjugate  Banach  spaces,  Studia  Math.  26
(1966), 307--313.

[G]  B.V.Godun, On weak$^{*}$ derived sets of sets of linear functionals,
Mat. Zametki 23 (1978), 607--616 (in Russian).

[GR] B.V.Godun and S.A.Rakov, Banach - Saks property and the  three
space problem,  Mat. Zametki. 31 (1982), 61-74. (in Russian).

[Gu] V.I.Gurarii, On openings and inclinations of  subspaces  of  a
Banach space, Teor.  Funktsii, Funktsional. Anal. i Prilozhen.
1 (1965), 194--204. (in Russian).

[HW] R.Herman and R.Whitley, An  example  concerning  reflexivity, 
Studia Math. 28 (1967), 289--294.

[JR] W.B.Johnson and H.P.Rosenthal, On $w^{*}$-basic sequences and  their
applications to the study of Banach  spaces,  Studia  Math.  43
(1972), 77--92.

[LT] J.Lindenstrauss and L.Tzafriri, Classical Banach  spaces I.
Sequence spaces (Springer, Berlin, 1977).

[Ma]  S.Mazurkiewicz,  Sur la  derivee  faible  d'un  ensemble   de
fonctionelles lineaires, Studia Math. 2 (1930), 68--71.

[Mc] O.C.McGehee, A proof of a statement of Banach about  the  weak$^{*}$
topology, Michigan Math. J. 15 (1968), 135--140.

[MM1] G.Metafune and V.B.Moscatelli, Generalized prequojections  and
bounded maps, Results in Math. 15 (1989), 172--178.

[MM2] G.Metafune and V.B.Moscatelli, Quojections and prequojections,
in: Advances in  the  Theory  of  Frechet  Spaces,  T.Terzioglu
(ed.), Kluwer, Dordrecht, 1989, 235--254.

[M1] V.B.Moscatelli, On strongly non-norming subspaces, Note Mat.  7
(1987), 311--314.

[M2]   V.B.Moscatelli,    Strongly    nonnorming    subspaces    and
prequojections, Studia Math. 95 (1990), 249--254.

[O1]  M.I.Ostrovskii, $w^{*}$-derived  sets  of  transfinite  order   of
subspaces of dual Banach spaces, Dokl. Akad. Nauk Ukrain. SSR,
Ser. A (1987), no. 10, 9--12 (in Russian).

[O2] M.I.Ostrovskii, On total nonnorming subspaces of a conjugate  Banach
space, Teor. Funktsii, Funktsional. Anal. i Prilozhen. 53 (1990),
119--123. (in Russian). Engl.
transl: J. Soviet Math. 58 (1992), no. 6, 577--579.

[O3]  M.I.Ostrovskii,  Regularizability   of
superpositions of inverse  linear  operators,  Teor.  Funktsii,
Funktsional. Anal i Prilozhen. 55 (1991), 96--100 (in Russian). Engl.
transl.: J. Soviet Math. 59 (1992), no. 1, 652--655.

[Pe] A.Pelczynski, Banach spaces on  which  every  unconditionally
converging operator is weakly  compact, Bull. l'Acad. Polon.
Sc. Ser. math. astr. phys. 10 (1962), 641--648.

[P]  Yu.I.Petunin, Conjugate Banach spaces containing  subspaces  of
zero  characteristic,  Dokl.  Akad.  Nauk  SSSR,  154   (1964),
527--529, Engl. transl.: Soviet Math. Dokl. 5 (1964), 131--133.

[PP] Yu.I.Petunin and A.N.Plichko, The theory of  characteristic  of
subspaces and its applications (Vyshcha Shkola, Kiev, 1980,  in
Russian).

[Pl] A.N.Plichko, On bounded biorthogonal systems in  some  function
spaces, Studia Math. 84 (1986), 25--37.

[S1] D.Sarason, On the order of a simply connected domain,  Michigan
Math. J. 15 (1968), 129--133.

[S2] D.Sarason, A remark on the weak-star  topology  of $l^{\infty }$,  Studia
Math. 30 (1968), 355--359.

[Sc] J.J.Schaffer,  Linear  differential  equations  and  functional
analysis. VI, Math. Ann. 145 (1962), 354--400.

[W]  P.Wojtaszczyk, Banach spaces for analysts, Cambridge  studies
in advanced mathematics 25 (Cambridge University Press, 1991).

\obeylines
Address: Institute for Low Temperature 
Physics and Engineering
Ukrainian Academy of Sciences
47 Lenin Avenue
Kharkov 310164 UKRAINE
\end{large}
\end{document}